\documentclass[final]{elsart}
\usepackage{amsmath,amssymb}
\usepackage{epsfig}

\newcommand{\bee}{\begin{equation}}
\newcommand{\eee}{\end{equation}}
\newcommand{\beq}{\begin{eqnarray*}}
\newcommand{\eeq}{\end{eqnarray*}}
\newcommand{\beql}{\begin{eqnarray}}
\newcommand{\eeql}{\end{eqnarray}}
\newcommand{\OO}{\mathcal{O}}
\newcommand{\R}{\mathbb R}
\newcommand{\C}{\mathbb C}

\newcommand{\eps}{\varepsilon}

\renewcommand{\>}{\rangle}

\begin{document}

\begin{frontmatter}
\title{Switching to nonhyperbolic cycles from codim 2 bifurcations of equilibria in ODEs}

\author[UU]{Yu.A. Kuznetsov}
\author[UT]{H.G.E. Meijer\corauthref{cor}}
\corauth[cor]{Corresponding author.}
\author[UG]{W. Govaerts}
\author[UG]{B. Sautois}

\address[UU]{Department of Mathematics, Utrecht University,Budapestlaan 6,
3584CD Utrecht,The Netherlands}
\address[UT]{Department of Applied Mathematics, Twente University, P.O. Box 217,
7500 AE, Enschede, The Netherlands}
\address[UG]{Department of Applied Mathematics and Computer Science, Ghent
University, Krijgslaan 281-S9, B-9000 Ghent, Belgium}

\begin{abstract}
The paper provides full algorithmic details on switching to the continuation
of all possible codim 1 cycle bifurcations from generic codim 2 equilibrium
bifurcation points in $n$-dimensional ODEs. We discuss the implementation
and the performance of the algorithm in several examples, including an
extended Lorenz-84 model and a laser system. 
\end{abstract}

\end{frontmatter}

\section{Introduction}
Consider a system of differential equations depending on two parameters
\begin{equation}
\dot{x}=f(x,\alpha),\ \ (x,\alpha) \in \R^n \times \R^2,
\label{eq:ODE}
\end{equation}
where $f$ is smooth.
In general, there are bifurcation curves in the $\alpha$-plane, at which 
the system exhibits codim 1 bifurcations, for example, fold or Hopf bifurcations of 
equilibrium points. Moreover, generically, one expects points of codim 2 bifurcations, 
where several curves corresponding to codim 1 bifurcations intersect transversally or 
tangentially. A codim 2 point is of particular interest if it is not only the origin
of some equilibrium bifurcation curves but also of some curves corresponding to bifurcations of
periodic orbits (cycles). Such points can be detected by purely local analysis of equilibria and
then be used to establish the existence of limit cycle bifurcations and other global phenomena 
that could hardly be proved otherwise. That is why codim 2 points are often called
the ``organizing centers" in applied literature.

The theory of codim 2 bifurcations of equilibria in generic systems ({\ref{eq:ODE})
is well-developed (see, for example,  \cite{Ar:83}, \cite{GuHo:83}, \cite{Ku:04}). 
There are five well-known codim 2 equilibrium bifurcations: cusp ({\sf CP}),
Bautin (generalized Hopf, {\sf GH}), double zero (Bodanov-Takens, {\sf BT}), zero-Hopf ({\sf ZH}), and
double Hopf ({\sf HH}). It follows from their analysis that branches of nonhyperbolic 
limit cycles can emanate from {\sf GH}, {\sf ZH}, and {\sf HH} points only. More precisely, 
a codim 1 bifurcation curve {\sf LPC}, along which a cycle with a nontrivial multiplier 
$\mu_1=1$ exists, emanates from a generic {\sf GH} point, while codim 1 bifurcation curves {\sf NS},
along which a cycle with a pair of multipliers $\mu_{1,2}= {\rm e}^{\pm i\theta}$
exists, are rooted at generic {\sf ZH} and {\sf HH} points. Notice that {\sf NS}
is used to denote both Neimark-Sacker and neutral saddle cycles where $\mu_{1}\mu_{2}=1$
and that no period-doubling curves can emanate from generic codim 2 equilibrium bifurcations.

Obviously, the application of these theoretical results to realistic models (\ref{eq:ODE})
is impossible without numerical tools. The numerical analysis of a codim 2 equilibrium 
bifurcation includes:
\begin{itemize}
\item detection and location of the point in a branch of a codim 1 bifurcation;
\item computation of the coefficients of the normal form of the restriction
of (\ref{eq:ODE}) to the critical center manifold at the bifurcation parameter values
and checking the nondegeneracy conditions;
\item verification of the transversality of the given family (\ref{eq:ODE}) to the
codim 2 bifurcation manifold and establishing a correspondence between the unfolding
parameters of the normal form and original system parameters $\alpha$;
\item computing accurate approximations of the codim 1 curves in the $\alpha$-space
and the corresponding singular orbits in the $x$-space near the bifurcation, sufficient 
to initialize the numerical continuation of these codim 1 curves using only
local information available at the codim 2 point.
\end{itemize}
While the first two problems were studied in detail (see, \cite{HANDBOOK} and references therein)
and have been implemented into the standard bifurcation software {\sc content} \cite{CONTENT}
and {\sc matcont} \cite{MATCONT}, two last issues received much less attention in the
numerical analysis literature, even if bifurcations of nonhyperbolic cycles are 
concerned. The present paper is aimed at bridging this gap by providing full algorithmic
details on switching to all possible codim 1 cycle bifurcations from generic {\sf GH}, {\sf ZH}, and 
{\sf HH} codim 2 points. 
 
One way to set up a computational switching procedure is to consider a 
smooth normal form for the codim 2 bifurcation including the parameters
$\beta \in \R^2$ 
\begin{equation}
\label{eq:nf}
\dot{w} = G(w,\beta), \qquad G:\R^{n_c} \times \R^2 \to \R^{n_c}.
\end{equation}
For all codim 2 equilibrium bifurcations these normal forms are known. 
Suppose that an exact or approximate formula is available that gives
the emanating codim 1 bifurcations for the normal form 
(\ref{eq:nf}). In order to transfer this to the original equation
(\ref{eq:ODE}) we need  a relation 
\begin{equation}
\label{eq:pars}
\alpha=V(\beta), \qquad V: \R^2 \to \R^2 
\end{equation} 
between the unfolding parameters $\beta$ and the given parameters
$\alpha$. In our context, V will be linearly approximated.
Moreover, we need a center manifold parametrization 
\begin{equation}
\label{eq:CM}
x=H(w,\beta), \qquad H: \R^{n_c} \times \R^2 \to \R^{n},
\end{equation}
that incorporates $\beta$. Taking (\ref{eq:pars}) and (\ref{eq:CM}) together as
$(x,\alpha)=(H(w,\beta),V(\beta))$ yields a center manifold for
the suspended system $\dot{x}=f(x,\alpha),\dot{\alpha}=0$.
The invariance condition for the center manifold now turns into
a {\em homological equation}:
\begin{equation}
\label{eq:HOM}
H_{w}(w,\beta)G(w,\beta) = f(H(w,\beta), V(\beta)),
\end{equation}
which we can solve by a recursive procedure based on Fredholm's solvability condition
that will give the Taylor coefficients of $G$ and $H$ with respect to $w$ and $\beta$.
We assume the Taylor series of $G$ to be known as
\[
G(w,\beta)=\sum_{|\nu|+|\mu| \geq 1} \frac{1} {\nu ! \mu !}
           g_{\nu \mu} w^{\nu} \beta^{\mu}, \quad
\]
and the Taylor series of $H$ and $V$ to be unknown
\[
H(w,\beta)=\sum_{|\nu|+|\mu| \geq 1} \frac{1} {\nu ! \mu !}
           h_{\nu \mu} w^{\nu} \beta^{\mu}, \quad
V(\beta)=\sum_{|\mu| \geq 1} \frac{1} {\mu !}
           v_{\mu} \beta^{\mu}.
\]
Here $\nu$ and $\mu$ are multi-indices. For  $\mu=0$ this 
reproduces the critical normal form coefficients first computed in \cite{Ku:99}, 
while the coefficients with $|\mu|\geq 1$ yield the necessary data on the 
parameter dependence. 

To summarize, a bifurcation point is detected within a certain small tolerance. As the
prediction depends on the initial point, this translates into small errors of the predicted
curve. If we start close enough to the actual new curve, any point will converge to it
and in general one expects a convergence cone \cite{JeDe:86}. If we parametrize the predicted
curve by $\eps$, the initial amplitude $\eps$ is to be chosen to be within the convergence cone,
see also Figure \ref{predictor}.
\begin{figure}\begin{center}
\includegraphics[width=8cm]{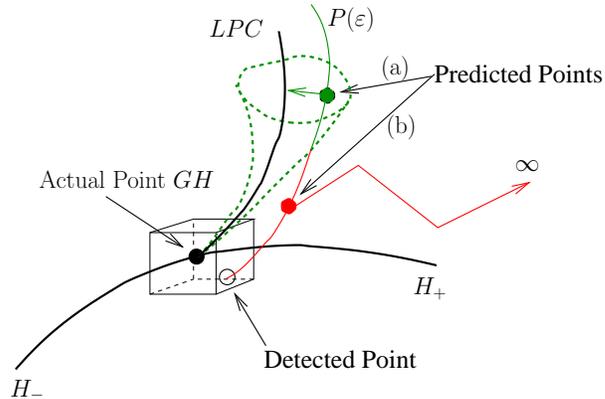}
\caption{Sketch of the switch in the case of a {\sf GH} bifurcation.
A predicted point along $P(\eps)$ (a) in the cone will converge to the LPC-curve,
outside (b) it will not.}
\label{predictor}
\end{center}\end{figure}

This procedure is adopted from \cite{HANDBOOK}, where it has been applied to the derivation
of the asymptotics of the fold and Hopf curves rooted at {\sf CP} and {\sf BT}
codim 2 points, as well as that for a homoclinic orbit to a saddle emanating
at the {\sf BT}-point. Recently, this technique has been
successfully used for switching at codim 2 fixed points of maps to the
continuation of nonhyperbolic periodic orbits rooted there \cite{GKKM:07}. 
Note that a similar procedure was suggested 
in \cite{IHS:98}, without using the Fredholm condition,
and carried through in the {\sf ZH}-case in \cite{IHS:00}, where, however, no 
asymptotics of codim 1 curves were derived. Finally, we point out that the 
problem of switching to the {\sf LPC}-curve at the {\sf GH} bifurcation has been briefly 
discussed in \cite{HANDBOOK} in a setting without the Taylor expansion in $\beta$. 

The paper is organized as follows. In Section \ref{Sec:AsymptCM} we revise
smooth parameter-dependent normal forms on center manifolds for the considered
codim 2 bifurcations, i.e. give $G(w,\beta)$ in {\sf GH}, {\sf ZH}, and 
{\sf HH} cases, and give the asymptotic expressions of the branches of nonhyperbolic
cycles in these normal forms. Then we perform the described above reduction procedure 
and derive the necessary coefficients $g_{\nu \mu}, h_{\nu \mu}$, and $v_{\mu}$ in terms
of $F$ and its derivatives. 
These coefficients are finally used to set up predictors for these branches in the
original system (\ref{eq:ODE}). An implementation of the resulting formulas in the
software {\sc matcont} is discussed at the end of this section. Section \ref{Sec:Examples}
presents several applications of the developed technique to known ODE models,
an extension of the Lorenz-84 system and a laser model, where we compare the asymptotic formulas for
the cycle bifurcations with numerically computed {\sf LPC}- and {\sf NS}-branches.
A discussion of existing results and open problems in switching to homoclinic
branches at {\sf ZH} and {\sf HH} bifurcations is given in Section \ref{Sec:Discussion}.

\section{Asymptotics and the Center Manifold}
\label{Sec:AsymptCM}

\subsection{The `new' curves}
The parameter-dependent normal forms are known and can be found in the standard
texts, e.g. \cite{Ku:04}. As the normal form and the asymptotic expressions
are the necessary theoretical ingredient, we present these here.

\subsubsection{Generalized Hopf}
Near a {\sf GH} bifurcation the vector field restricted to the
center manifold is given by
\begin{equation}\label{nf-GH}
\dot{w} = \lambda(\beta)w + c_{1}(\beta)w|w|^{2} + c_{2}(\beta_{2})w|w|^{4}+\OO(|w|^{6}),~~ w \in
\C,
\end{equation}
where $\lambda(0)=i\omega$, and this bifurcation is characterized by 
$d_{1}=\Re(c_{1}(0))=0$ and $d_{2}=\Re(c_{2}(0))\neq 0$. A curve {\sf LPC} of fold
bifurcation of limit cycles emanates from this point. Let us write
$w=\rho {\rm e}^{i\psi}$, $\lambda(\beta) = i\omega + \beta_{1}+ib_{1}(\beta)+O(|\beta|^{2})$
and $\Re(c_{1}(\beta)) = \beta_{2}+O(|\beta|^{2})$. If we now truncate the
normal form to fifth order in $w$, then the curve {\sf LPC}
is given by 
\begin{equation}\label{asymp-GH}
\rho = \eps, \beta_{1} = d_{2}\eps^{4}, \beta_{2} = -2d_{2}\eps^{2}.
\end{equation}

\subsubsection{Zero-Hopf}
Near a {\sf ZH} bifurcation the vector field restricted to the
center manifold is given by
\begin{equation}\label{nf-ZH}
\left(\begin{array}{c}\dot{x} \\ \dot{w} \end{array}\right) =
\left(\begin{array}{c} 
\beta_{1} + f_{200}x^{2} + f_{011}|w|^{2} + f_{300}x^{3} +f_{111}x|w|^{2} \\
(i\omega(\beta)+\beta_{2})w + g_{110}xw + g_{210}x^{2}w + g_{021}w|w|^{2} 
\end{array}\right)+\OO(\|(x,w)\|^{4}),
\end{equation}
where $(x,w) \in \R \times \C$. An extra Neimark-Sacker (torus) bifurcation of limit cycles 
({\sf NS}) occurs if $\Re(g_{110})f_{011}<0$.

The asymptotic expression is 
\begin{equation}\label{asymp-ZH:a}\begin{array}{c}
{\displaystyle \rho = \eps, x = -\frac{f_{111}+2g_{021}}{2f_{200}}\eps^{2}, 
\beta_{1} = -f_{011}\eps^{2}},\\
{\displaystyle \beta_{2} = \frac{2(\Re(g_{110})-f_{200})\Re(g_{021})
+\Re(g_{110})f_{111}}{2f_{200}}\eps^{2}}.
\end{array}\end{equation}
This agrees with a formula given in \cite{Ga:93}.

\subsubsection{Double-Hopf}
For a {\sf HH} bifurcation the dynamics on the center manifold is
governed by the following normal form:
\begin{equation}\label{nf-HH}
\left(\begin{array}{c}\dot{w_{1}} \\ \dot{w_{2}} \end{array}\right) =
\left(\begin{array}{l} 
(i\omega_{1}(\beta)+\beta_{1})w_{1} + f_{2100}w_{1}|w_{1}|^{2} + f_{1011}w_{1}|w_{2}|^{2}\\ 
(i\omega_{2}(\beta)+\beta_{2})w_{2} + g_{1110}w_{2}|w_{1}|^{2} + g_{0021}w_{2}|w_{2}|^{2} 
\end{array}\right)+\OO(\|(w_{1},w_{2})\|^{4}),
\end{equation}
where $(w_1,w_2) \in \C \times \C$.
Then there are generically two half-lines along which there is a
{\sf NS} bifurcation of limit cycles. In polar coordinates
$w_{1} = \rho_{1}{\rm e}^{i\psi_{1}},w_{2} = \rho_{2}{\rm e}^{i\psi_{2}}$
their asymptotics are given as
\begin{eqnarray}\label{asymp-HH} 
(\rho_{1},\rho_{2},\beta_{1},\beta_{2}) &=& \left(\eps,0,-\Re(f_{2100})\eps^2,
-\Re(g_{1110})\eps^2\right),\\ 
(\rho_{1},\rho_{2},\beta_{1},\beta_{2}) &=& \left(0,\eps,-\Re(f_{1011})\eps^2,
-\Re(g_{0021})\eps^2\right).
\end{eqnarray}

\subsection{Coefficients of parameter-dependent center manifolds}
We assume that the critical normal form coefficients are known 
(see \cite{Ku:99} and \cite{HANDBOOK}) and give here only parameter-related 
coefficients $h_{\nu \mu}$ from the homological equation. These provide
in each case a linear approximation to the parameter transformation (\ref{eq:pars}) .

\subsubsection{Generalized Hopf}
Here we closely follow the idea outlined in \cite{HANDBOOK}. We first expand the
eigenvalue and the first Lyapunov coefficient in the original parameters $\alpha$ 
and collect the equations to obtain the transformation to the unfolding parameters $\beta$.
Alternatively, one can normalize already in (\ref{GH:eqs1}) to obtain an
orthogonal frame from these equations and obtain scalings from the higher
order equations (\ref{GH:eqs2}). Below we have $\mu=(10),(01)$ as indices
and $v_{10}=(1,0)$, $v_{01}=(0,1)$ as vectors.

The first two equations (actually four) coming from (\ref{eq:HOM}) are
\begin{equation}\begin{array}{rcl}\label{GH:eqs1}
Ah_{00\mu} &=& -J_{1} v_{\mu},\\
(A-i\omega I_{n}) h_{10\mu} &=& \gamma_{1,\mu}q - A_{1}(q,v_{\mu}) - B(q,h_{00\mu})
\end{array}\end{equation}
The first equation is nonsingular and from the second we find $\gamma_{1,\mu}$ using the
Fredholm alternative. The other systems from (\ref{eq:HOM}) are
\begin{equation}\begin{array}{rcl}\label{GH:eqs2}
(A-2i\omega I_{n})h_{20\mu} &=& 2h_{2000}\gamma_{1,\mu} -\left[C(q,q,h_{00\mu})+2B(q,h_{10\mu})+B(h_{2000},h_{00\mu})\right. \\
 & & \left. + B_{1}(q,q,v_{\mu})+A_{1}(h_{2000},v_{\mu})\right] ,\\
Ah_{11\mu} &=& 2\Re(\gamma_{1,\mu})h_{1100} -\left[ C(q,\bar{q},h_{00\mu})+B(h_{1100},h_{00\mu})\right. \\
 & &\left. +B(\bar{q},h_{10\mu})+B(q,h_{01\mu})+ B_{1}(q,\bar{q},v_{\mu})+A_{1}(h_{1100},v_{\mu})\right],\\
(A-i\omega I_{n})h_{21\mu} &=& 2\gamma_{2,\mu}q +h_{2100}(2\gamma_{1,\mu}+\bar{\gamma}_{1,\mu})+2h_{10\mu}c_{1}\\
& & -\left[D(q,q,\bar{q},h_{00\mu})+2C(q,h_{1100},h_{00\mu})+2C(q,\bar{q},h_{10\mu})\right. \\
 & &+C(q,q,h_{01\mu})+C(h_{2000},\bar{q},h_{00\mu})+ 2B(q,h_{11\mu})\\
 & &+2B(h_{1100},h_{10\mu})+B(h_{2000},h_{01\mu})+B(h_{2100},h_{00\mu})\\
 & &+B(h_{20\mu},\bar{q})+C_{1}(q,q,\bar{q},v_{\mu})+2B_{1}(h_{1100},q,v_{\mu})\\
 & &\left. +B_{1}(h_{2000},\bar{q},v_{\mu})+A_{1}(h_{2100},v_{\mu})\right],\\
\end{array}\end{equation}
The first two are nonsingular and with the Fredholm alternative we find $\gamma_{2,\mu}$. The
parameter transformation (\ref{eq:pars}) is given by
\begin{equation}\label{eq:pars-GH}
\alpha=\left(\Re\left(\begin{array}{cc}\gamma_{1,10} & \gamma_{1,01} \\
\gamma_{2,10} & \gamma_{2,01} \end{array}\right)\right)^{-1}\beta.
\end{equation}

\subsubsection{Zero-Hopf}
This case is also treated in \cite{IHS:00}, however with only one parameter and
for hyperbolic periodic orbits. Thus our computational scheme is different.
We list only the necessary equations. 
\begin{equation}\label{cmred-ZH}\begin{array}{lrcl}
{\rm (a)} &A[h_{00010},\ h_{00001}] &=& [q_{1},\ 0]-J_{1} [v_{10},\ v_{01}],\\
{\rm (b)} &A[h_{10010},\ h_{10001}] &=& [h_{20000},\ 0] -A_{1}(q_{1},[v_{10},\ v_{01}]) \\ & & &- B(q_{1},[h_{00010},\ h_{00001}])\\
{\rm (c)} &(A-i\omega I_{n}) [h_{01010},\ h_{01001}] &=& [h_{11000},\ q_{2}] - A_{1}(q_{2},[v_{10},\ v_{01}]) \\
 & & &- B(q_{2},[h_{00010},\ h_{00001}])
\end{array}\end{equation}
In contrast to the other cases, here the first system is already singular. Taking the inner-product with
the adjoint null-vector we obtain the new orthogonal frame
\begin{equation}\begin{array}{c}
\gamma=(\gamma_{1},\gamma_{2}) = p_{1}^{T}J_{1},\quad s_{1}^{T}=\gamma/\|\gamma\|^{2},s_{2}^{T} = (-\gamma_{2},\gamma_{1}),\\
v_{10} = s_{1} +\delta_{1}s_{2},v_{01} = \delta_{2}s_{2}.
\end{array}\end{equation}
Polynomial terms in the normal form (\ref{nf-ZH}) like $\beta_{1}x$ are also resonant, but they can be eliminated by
hypernormalization. After solving (\ref{cmred-ZH}.a) with a bordered matrix, see \cite{Go:00}, still a multiple
of $q_{1}$ may be added to $h_{00010}$. We use this to perform hypernormalization.
Let us write
$$
r_{1} = -A^{INV}\left(\begin{array}{c} q_{1}-J_{1}s_{1}\\ 0\end{array}\right),\quad 
r_{2} = -A^{INV}\left(\begin{array}{c} -J_{1}s_{2}\\ 0\end{array}\right),
$$
where $A^{INV}$ indicates the use of the bordered matrix, then we can write
$$
h_{00010} = r_{1} + \delta_{1}r_{2}+\delta_{3}q_{1}, \quad 
h_{00001} =\delta_{2}r_{2}+\delta_{4}q_{1},
$$
for some $\delta$'s.
Then by applying the Fredholm alternative to (\ref{cmred-ZH}.b,c) we can solve for all $\delta$'s at once.
\begin{equation}\begin{array}{rcl}
LL \left(\begin{array}{c}\delta_{1} \\ \delta_{3} \end{array}\right) &=& -\left(\begin{array}{c} 
\<p_{1},A_{1}(q_{1},r_{1}) + B(q_{1},r_{1})\> \\ 
\<p_{2},A_{1}(q_{2},r_{1}) + B(q_{2},r_{1})\> \end{array}\right) \\
\Re(LL) \left(\begin{array}{c}\delta_{2} \\ \delta_{4} \end{array}\right) &=& \left(\begin{array}{c} 0 \\ 1 \end{array}\right)
\end{array}\end{equation}
where
$$
LL= \left(\begin{array}{cc}
\<p_{1},A_{1}(q_{1},r_{2}) + B(q_{1},r_{2})\> & 2f_{200}\\
\<p_{2},A_{1}(q_{2},r_{2}) + B(q_{2},r_{2})\> & g_{110}
\end{array}\right).
$$

\subsubsection{Double Hopf}
Although high-dimensional, this case can be treated in a relatively simple manner.
Using the same notation as for the generalized Hopf from (\ref{eq:HOM}) we get 
\begin{equation}\begin{array}{rcl}
Ah_{0000\mu} &=& -J_{1} v_{\mu},\\
(A-i\omega_{1} I_{n}) h_{1000\mu} &=& \gamma_{1,\mu}q_{1} - A_{1}(q_{1},v_{\mu}) - B(q_{1},h_{0000\mu}),\\
(A-i\omega_{2} I_{n}) h_{0010\mu} &=& \gamma_{2,\mu}q_{2} - A_{1}(q_{2},v_{\mu}) - B(q_{2},h_{0000\mu}).
\end{array}\end{equation}
As the first equation is non-singular, formal substitution of $h_{000010}$ and
$h_{000001}$ and the Fredholm alternative leads to the same transformation (\ref{eq:pars-GH})
from unfolding to the system parameters.

\subsection{Implementation of the Predictors}

We have implemented our switching routines in {\sc matcont} \cite{MATCONT}. For the continuation of 
{\sf LPC} and {\sf NS} curves it uses a minimally augmented defining system \cite{KuGoDoDh:2005}, i.e.
we need to supply an approximation of the limit cycle, the period and the parameters. 
The parameters follow from applying the inverse transformation to (\ref{eq:pars}).
There is always one dynamic variable $\psi$ giving a free phase shift along the bifurcating
limit cycle with a period $\frac{2\pi}{\omega_{1}(\eps)}$. For the initial cycle we
make an equidistant mesh $\psi=2n\pi/N, n=0\dots N$ where $N+1$ is the number of mesh
points. Let $q$ denote the eigenvector corresponding to the eigenvalue $i\omega_{1}$,
then points on the limit cycle are given by $x_{0}+\eps(qe^{i\psi}+\bar{q}e^{-i\psi})$.
Similarly, terms as $\eps^{2}h_{20}e^{2i\psi}$ and $\eps^{2}h_{0010}$ are included.
An internal routine of {\sc matcont} then adapts this limit cycle on an equidistant
mesh to a mesh defined at the non-equidistant collocation points.

For the {\sf NS} curves the system is augmented with the real part $k$ of the multiplier.
In this case the normal forms (\ref{nf-ZH}),(\ref{nf-HH}) also define a second rotation
with frequency $\omega_{2}(\eps)$ and we have $k=\cos\left(\frac{2\pi\omega_{2}(\eps)}{\omega_{1}(\eps)}\right)$.

{\sc matcont} uses Moore-Penrose continuation for which also a tangent vector to the
bifurcation curve is needed. This tangent vector is easily obtained by differentiating
the predictor w.r.t. $\eps$. 

Below we list some case-specific details.

\subsubsection{Generalized Hopf}
The period is given by $T=2\pi/\omega+(2d_{2}b_{1,2}-\Im(c_{1}(0)))\eps^{2}$,
with $b_{1,2}=\frac{\partial b_{1}}{\partial\beta_{2}}$. The parameters are given by
$\alpha=\alpha_{0}+V(0,-2d_{2}\eps^{2})^{T}$.

Note that for a $\eps^{4}$-approximation also seventh order derivatives would be
needed; this follows from Remark 3.3.2 in \cite{Me:06}. Therefor we restrict 
to $\OO(\eps^{3})$ in the implementation.

\subsubsection{Zero-Hopf}
In the continuation we also need to provide the period and the multiplier. Approximating
formulas are defined as follows where $x,\beta_{1},\beta_{2}$ are as in (\ref{asymp-ZH:a})
\begin{equation}\label{asymp-ZH:b}\begin{array}{c}
T=2\pi/\omega(0)-(\omega_{1}\beta_{1}+\omega_{2}\beta_{2}+\Im(g_{110})x)-\Im(g_{021})\eps^{2},\\
k=1-(4\pi \Re(g_{110})f_{011})(\eps/\omega_{0})^{2}.
\end{array}\end{equation}

\subsubsection{Double Hopf}
Approximating formulas for the period and the multiplier on one branch are given by
\begin{equation}\label{asymp-HH:b}
\begin{array}{l}
T=\frac{2\pi}{\omega_{1}+d\omega_{1}\eps^{2}},\qquad k=\cos(T(\omega_{2}+d\omega_{2}\eps^{2})),\\
(d\omega_{1},d\omega_{2}) = -\Im(\gamma_{1} \gamma_{2})^{T}(\Re(\gamma_{1} \gamma_{2})^{T})^{-1}
\Re(f_{2100},g_{1110})^{T}+\Im(f_{2100},g_{1110}).
\end{array}
\end{equation}
and similarly for the other branch.

\section{Examples}
\label{Sec:Examples}

\subsection{New curves in an extension of the  Lorenz-84 model}
\label{Sec:Lorenz-84}
The first example is an extended version of the Lorenz-84 model. A bifurcation analysis
of this model was presented in \cite{ShNiNi:95,Ve:03}. In this system $X$ models the
intensity of a baroclinic wave and $Y$ and $Z$ the sine and cosine coefficients of the wave.
This model may be extended with a variable $U$ to study the influence of external
parameters such as temperature and the model then shows several limit cycle
bifurcations \cite{KuMeVe:04}. It has the form:
\begin{equation}
\left\{\begin{array}{rcl}
\dot{X} & = & -Y^{2}-Z^{2}-\alpha X + \alpha F - \gamma U^2  \\
\dot{Y} & = & X Y- \beta X Z - Y + G \\
\dot{Z} & = & \beta X Y + X Z - Z \\
\dot{U} & = & -\delta U +\gamma UX + T
\end{array}
\right.
\label{extLor}
\end{equation}
The parameters $F$ and $T$ are varied while we fix $\alpha=.25,\beta=1,G=.25,\delta=1.04,\gamma=.987$.
The bifurcation diagram displays one fold bifurcation and two Hopf bifurcation curves, see Figure
\ref{lorenzdiagram}. We find all codim 2 points of equilibria, in particular $GH,ZH$ and
$HH$.


We have applied our switching routines to all three emanating curves, since 
the $NS$ bifurcation from ZH is a neutral saddle. The predictions
in parameterspace are shown in Figure \ref{lorenzdiagram} next to the numerically continued curves. The
predicted points were used as a starting point point for the continuation of these limit cycle
bifurcations, which shows that our approach works. Another numerical check is provided by inspecting
the tangent vector, which we provide together with a first point. When we find a second point
on the curve by continuation and adapt the defining system, we will obtain a more precise tangent vector.
For a small continuation step, this tangent vector and the predicted one should be close.
Indeed, for the examples reported here, the first digits always coincided.

\begin{table}[t]\begin{center}
\begin{tabular}{|c|c|c|l|}
\hline \hline
Label & F & T & Normal Form coefficients \\
\hline
GH & $2.3763601$ & $.050197432$ & $d_{2}=0.1558012$\\
HH & $2.5332211$ & $.026273943$ & $p_{11}p_{22}=-1$, $\theta=-3.648550$, $\delta=-1.052987$\\
& & & $\Theta=1230.630$, $\Delta=-210.861$\\
ZH & $1.2834193$ & $.000126541$ & $s=1$, $\theta=0.3715145$, $E=-1$\\
\hline \hline
\end{tabular}
\end{center}
\caption{Parameter values of $F$ and $T$ at the bifurcation points in Figure 
\ref{lorenzdiagram} together with normal coefficients (scaled, see \cite{Ku:04}).}
\end{table}

\begin{figure}
\begin{center}
\includegraphics[width=13.8cm]{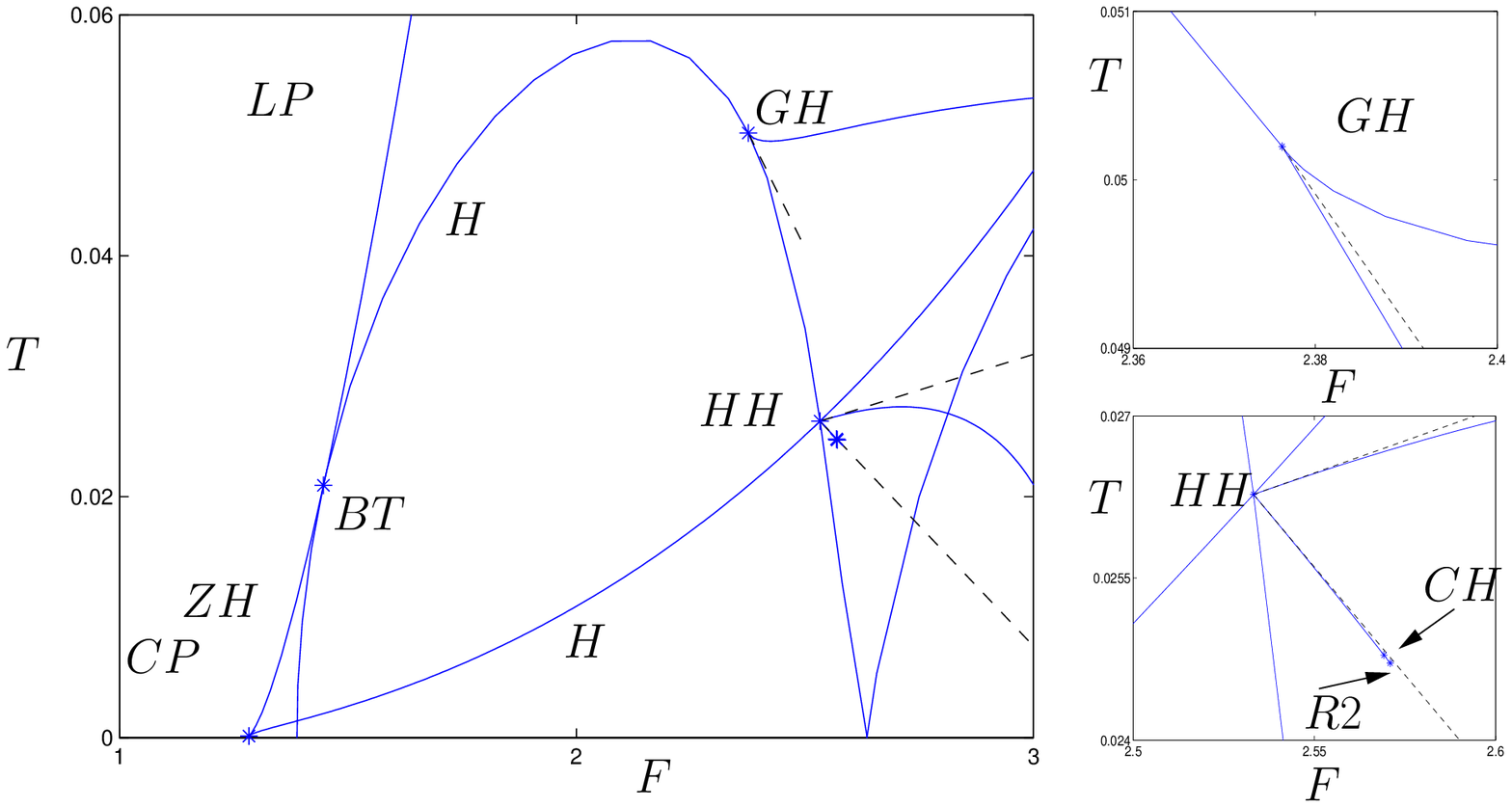}
\caption{Bifurcation diagram of the Extended Lorenz-84 model. Symbols denote $LP$ Limit Point,
$H$ Hopf, $LPC$ Limit Point of Cycles, $NS$ Neimark-Sacker, $GH$ Generalized Hopf, $HH$ Double Hopf,
$ZH$ Zero-Hopf, $BT$=Bogdanov-Takens. Dashed lines show the predicted new curves; (a) Zoom near the $GH$ point,
(b) Zoom near the $HH$ point.}
\label{lorenzdiagram}
\end{center}
\end{figure}

Finally we present some measure of the error of the switching routines as a function of the
initial amplitude $\eps$, see Figure \ref{lorenzerror} and its caption. Interestingly this Figure
represents the idea of Figure \ref{predictor}. Using a small initial amplitude $\eps$
may not work due to a numerical error in the calculated codim 2 point, on the other
hand $\eps$ must not be taken too large for the approximation to remain valid.
\begin{figure}[ht!]
\includegraphics[width=6cm]{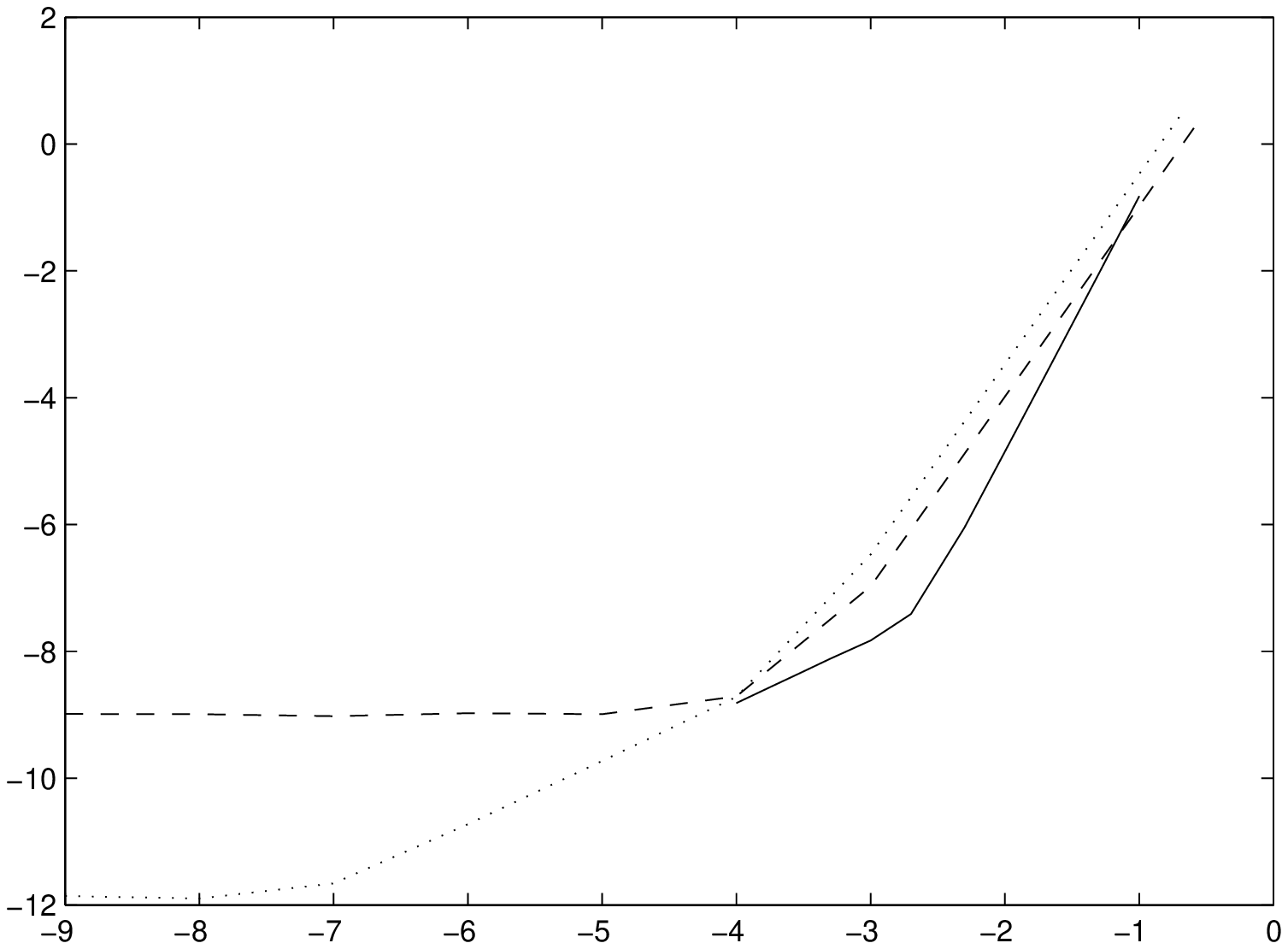}
\includegraphics[width=6cm]{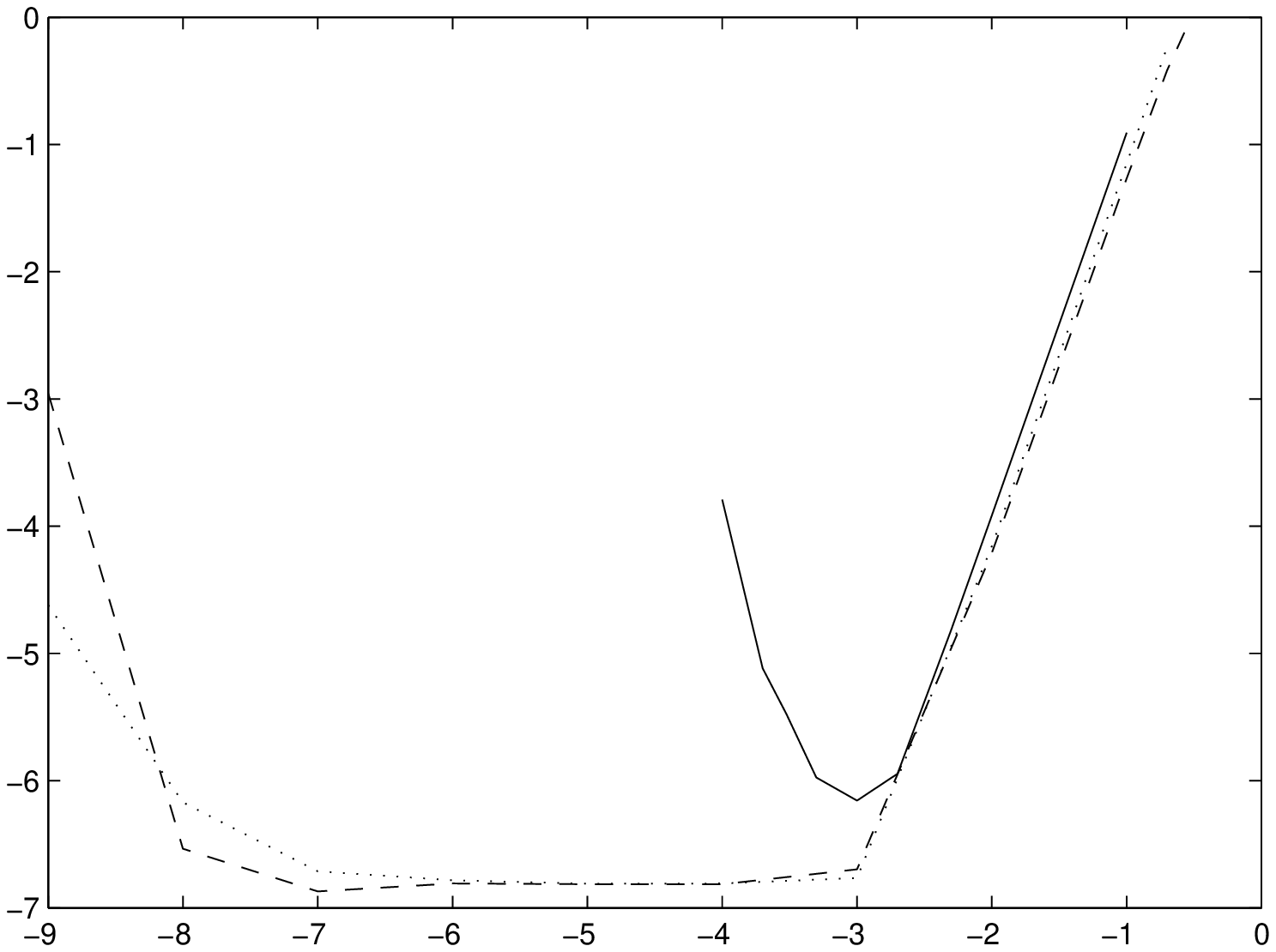}
\caption{Error measures:(a) The residual $R$ of the first Newton-step,
(b) The distance between the predicted and the first corrected point. For this we left
out the tangent vector in the continuation. {\sc matcont} then tries to correct the
first point immediately instead of starting the continuation. We have taken 20 mesh
and 4 collocation points and $\eps\in[10^{-7},.2]$. Data shown if predictor converged.}
\label{lorenzerror}
\end{figure}

\subsection{Switching in a Laser model}
\label{Sec:Laser}
In \cite{WiCh:06} a single-mode inversionless laser with a three-level phaser was studied and
shown to operate in various modes. These modes are ``off'' (non-lasing), continuous waves,
periodic, quasi-periodic and chaotic lasing. In particular the boundary of the region of chaos
seems to be defined by several limit cycle bifurcations born from several codim 2 equilibrium
bifurcations. Thus we want to start with our routines such boundary computations without first
doing simulations and limit cycle continuations in this 9-dimensional system.

The model is 9-dimensional system given by 3 real and 3 complex equations:
\begin{equation}\label{Laser}
\left\{\begin{array}{rcl}
\dot{\Omega}_{l}&=&-\frac{\gamma_{cav}}{2}\Omega_{l}-g\Im(\sigma_{ab}),\\
\dot{\rho}_{aa}&=&R_{a}-\frac{i}{2}(\Omega_{l}(\sigma_{ab}-\sigma^{*}_{ab})+\Omega_{p}(\sigma_{ac}-\sigma^{*}_{ac})),\\
\dot{\rho}_{bb}&=&R_{b}+\frac{i}{2}\Omega_{l}(\sigma_{ab}-\sigma^{*}_{ab}),\\
\dot{\sigma}_{ab}&=&-(\gamma_{1}+i\Delta_{l})\sigma_{ab}-\frac{i}{2}(\Omega_{l}(\rho_{aa}-\rho_{bb})-\Omega_{p}\sigma_{cb}),\\
\dot{\sigma}_{ac}&=&-(\gamma_{2}+i\Delta_{p})\sigma_{ac}-\frac{i}{2}(\Omega_{p}(2\rho_{aa}+\rho_{bb}-1)-\Omega_{l}\sigma^{*}_{cb}),\\
\dot{\sigma}_{cb}&=&-(\gamma_{3}+i(\Delta_{l}-\Delta_{p}))\sigma_{cb}-\frac{i}{2}(\Omega_{l}\sigma^{*}_{ac}-\Omega_{p}\sigma_{ab}),\\
\end{array}
\right.
\end{equation}
with $R_{a}=-.505\rho_{aa}-.405\rho_{bb}+.45,R_{b}=.0495\rho_{aa}-.0505\rho_{bb}+.0055$ 
and $\Delta_{l}:=\Delta_{cav}+g\Re(\sigma_{ab})/\Omega_{l}$. The parameters are fixed at
$\gamma_{1}=.05,\gamma_{2}=.25525,\gamma_{3}=.25025,\gamma_{cav}=.03,g=100,\Delta_{p}=0$
while $\Delta_{cav}$ and $\Omega_{p}$ are varied to study several detuning effects.
For more details, see \cite{WiCh:06}. 

\begin{figure}[ht!]
\includegraphics[width=12cm]{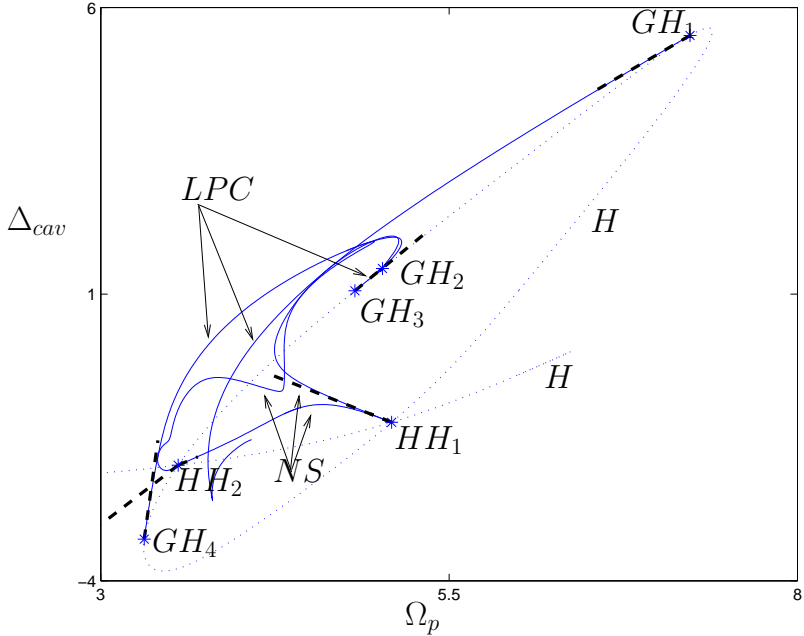}
\caption{Bifurcation diagram of the inversionless laser. Hopf curves(denoted by $H$) are dotted,
Limit cycle bifurcations are denoted by $LPC$ Limit Point of Cycles, $NS$ Neimark-Sacker.
Dashed lines show the predicted new curves.}
\label{laserdiagram}
\end{figure}

We have reproduced a part from the bifurcation diagram
which corresponds to continuous wave and periodically pulsating solutions, i.e. with
$\Omega_{l}\neq 0$, see Figure \ref{laserdiagram}. As the system has $\mathbb{Z}_{2}$-symmetry
the same bifurcations are found for $\Delta_{cav}\rightarrow-\Delta_{cav}$. For clarity
of the figure we do not display these here. We list the codim 2 points in Table \ref{codim2laser}.
The normal form coefficients of HH$_{1}$ confirm the claim of \cite{WiCh:06} that the most complicated type was
encountered; only the 3-torus is (un)stable. This is also confirmed
when we continue the Neimark-Sacker bifurcations. For HH$_{2}$ the NS curves are not
in the same quadrant defined by the Hopf curves, while they are for HH$_{1}$.
All cycle bifurcations where computed with 20 mesh points and 4 collocation points and
the initial amplitude was set to $\eps=.001$, which worked immediately in all cases.
Let us remark that one LPC curve connects GH$_{2}$ and $GH_{3}$ points and stays
close to the Hopf curve. Similarly, a NS curve starts at HH$_{1}$, becomes neutral between
two 1:2 resonances and ends at HH$_{2}$. It would have taken much more effort to
find this feature otherwise.

\begin{table}[t]\begin{center}\label{codim2laser}
\begin{tabular}{|c|c|c|l|}
\hline \hline
Label & $\Omega_{p}$ & $\Delta_{cav}$ & Normal Form coefficients \\
\hline
$GH_{1}$ & $7.228819$ & $5.511455$ & $d_{2}=-46.49852$\\
$GH_{2}$ & $5.021574$ & $1.446387$ & $d_{2}=3.813132$\\
$GH_{3}$ & $4.824066$ & $1.059367$ & $d_{2}=195.1119$\\
$GH_{4}$ & $3.312120$ & $-3.273568$ & $d_{2}=-6.468468$\\
$HH_{1}$ & $5.087299$ & $-1.2362053$ & $p_{11}p_{22}=-1$, $\theta=-.07194543$, $\delta=-13.91412$\\
& & & $\Theta=.9595389$, $\Delta=-2602.275$\\
$HH_{2}$ & $3.555848$ & $-1.983857$ & $p_{11}p_{22}=1$, $\theta=-.1179924$, $\delta=-26.59452$\\
& & & $\Theta=-10.81042$, $\Delta=-2713.608$\\
\hline \hline
\end{tabular}
\caption{Parameter values of $\Omega_{p}$ and $\Delta_{cav}$ at the codim 2 points together with normal coefficients 
(scaled, see \cite{Ku:04}).}
\end{center}\end{table}

\begin{figure}[ht!]
\includegraphics[width=13cm]{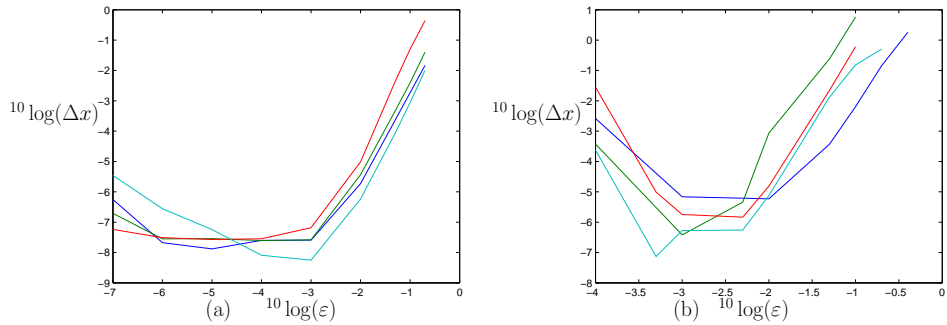}
\caption{Error measures along the eight curves: $^{10}\log$ of the
distance between the predicted and the first corrected point versus
$^{10}\log(\eps)$. (a) Along the Neimark-Sacker curves (b) Along the
LPC curves. This Figure again resembles the idea of Figure \ref{predictor}.}
\label{lasererror}
\end{figure}

\section{Discussion}
\label{Sec:Discussion}
This paper contributes to the bifurcation analysis of codim 2 singularities 
of equilibria in multidimensional ODEs by providing explicit predictors 
for branches of nonhyperbolic cycles emanating from these bifurcations. 
We have tested it on several examples with good results. We believe that
this work will further facilitate automated analysis of nonlinear systems.
However, we like to mention that we also tried the double Hopf point in a model for the 
lateral pyloric neuron \cite{GGK:97,GKS:00}. Although we were able to switch to
one branch and continue it without any problem, the Jacobian of the defining
system along the second branch was numerically singular. In this model with
multiple time scales probably a special numerical scheme is necessary.
 
It is well known that branches of orbits homoclinic to hyperbolic equilibria are
also rooted at {\sf BT}, {\sf ZH}, and {\sf HH} codim 2 bifurcation points.
The {\sf BT} case has been treated in \cite{Be:94} (see also \cite{HANDBOOK},
where the computational setting is most close to the present paper). The 
corresponding predictor for the homoclinic branch is implemented in {\sc matcont}.
The problem of providing predictors for homoclinic branches rooted at {\sf ZH} and {\sf HH}
points is more challenging. Some important results in this direction are
obtained in \cite{BrVe:84,Ga:93,ChKi:04}, where the systems reduced to the
center manifold were considered. However, a complete set of formulas
suitable for switching to homoclinic curves in these cases is still not available.
For instance, in the {\sf ZH} case the normal form (\ref{nf-ZH}) exhibits
homoclinic bifurcations of saddle-focus equilibria in the parameter plane
along a bifurcation curve with the linear approximation
\begin{eqnarray*}
&&\hspace{-1.0cm}\beta_{2,hom} = 
\frac{\Re(g_{110})\beta_{1}}{f_{200}(2f_{200}-3\Re(g_{110}))}\left[
\Re(g_{210}) - \frac{3\Re(g_{110})}{2f_{200}}f_{300} + 
\frac{(f_{200}-\Re(g_{110}))}{f_{011}}f_{111}\right.\\
&&~~~~~~~~~~~~~~~~~~~~~~~~~~~~~~~~~~~~~~~~~~~~~\left.
-\frac{2(f_{200}-\Re(g_{110}))^{2}}{f_{011}\Re(g_{110})}\Re(g_{021})\right],
\end{eqnarray*}
provided that $\Re(g_{110})f_{011}<0$ and $\Re(g_{110})f_{200}<0$.
Application of the above reduction to the parameter-dependent center manifolds
in the {\sf ZH} case yields an approximation to the bifurcation curve in the parameter plane.
Now the challenge is to construct a suitable initial solution in state space.
On this work in progress will be reported elsewhere.

Another direction for future research is a problem of switching to secondary cycle
bifurcations at codim 2 bifurcations of cycles in (\ref{eq:ODE}). Here a generalization of
the periodic normalization technique from \cite{KuGoDoDh:2005} to critical codim 2 cases and 
its extension to parameter-dependent systems in the spirit of \cite{GKKM:07} are required.

\section*{Acknowledgement}
The authors want to thank S. Wiezcorek for bringing up and his assistance with the laser model.

\bibliographystyle{plain}
\bibliography{cycle-switch}

\begin{thebibliography}{10}

\bibitem{Ar:83}
V.I. Arnold.
\newblock {\em {G}eometrical {M}ethods in the {T}heory of {O}rdinary
  {D}ifferential {E}quations}.
\newblock {S}pringer-{V}erlag, {N}ew {Y}ork, {H}eidelberg, {B}erlin, 1983.

\bibitem{Be:94}
{\mbox{W.-J}}.~Beyn.
\newblock Numerical analysis of homoclinic orbits emanating from a
  {T}akens-{B}ogdanov point.
\newblock {\em IMA J. Numer. Anal.}, 14:381--410, 1994.

\bibitem{HANDBOOK}
{\mbox{W.-J}}.~Beyn, A.~Champneys, E.~Doedel, W.~Govaerts,
  {\mbox{Yu.A}}.~Kuznetsov, and B.~Sandstede.
\newblock Numerical continuation, and computation of normal forms.
\newblock In B.~Fiedler, editor, {\em Handbook of {D}ynamical {S}ystems, {V}ol.
  2}, pages 149--219. Elsevier Science, Amsterdam, 2002.

\bibitem{BrVe:84}
H.~W. Broer and G.~Vegter.
\newblock Subordinate \v{S}il'nikov bifurcations near some singularities of
  vector fields having low codimension.
\newblock {\em Ergodic Theory Dynam. Systems}, 4:509--525, 1984.

\bibitem{ChKi:04}
A.~R. Champneys and V.~Kirk.
\newblock The entwined wiggling of homoclinic curves emerging from
  saddle-node/{H}opf instabilities.
\newblock {\em Phys. D}, 195:77--105, 2004.

\bibitem{MATCONT}
A.~Dhooge, W.~Govaerts, and {\mbox{Yu.A}}.~Kuznetsov.
\newblock {\sc matcont}:{A} {\sc matlab} package for numerical bifurcation
  analysis of {ODE}s.
\newblock {\em ACM Trans. Math. Software}, 29:141--164, 2003.

\bibitem{Ga:93}
P.~Gaspard.
\newblock Local birth of homoclinic chaos.
\newblock {\em Phys. D}, 62:94--122, 1993.

\bibitem{GKKM:07}
{\mbox{R.K}}.~Ghaziani, W.~Govaerts, {\mbox{Yu.A}}.~Kuznetsov, and
  {\mbox{H.G.E}}.~Meijer.
\newblock Numerical methods for two-parameter local bifurcation analysis of
  maps.
\newblock {\em SIAM J. Sci. Comput.}, 29:2644--2667, 2007.

\bibitem{GGK:97}
W.~Govaerts, J.~Guckenheimer, and A.~Khibnik.
\newblock Defining functions for multiple {H}opf bifurcations.
\newblock {\em SIAM J. Numer. Anal.}, 34(3):1269--1288, 1997.

\bibitem{GKS:00}
W.~Govaerts, Yu.~A. Kuznetsov, and B.~Sijnave.
\newblock Numerical methods for the generalized {H}opf bifurcation.
\newblock {\em SIAM J. Numer. Anal.}, 38(1):329--346, 2000.

\bibitem{Go:00}
{\mbox{W.J.F}}.~Govaerts.
\newblock {\em Numerical {M}ethods for {B}ifurcations of {D}ynamical
  {E}quilibria}.
\newblock SIAM, Philadelphia, 2000.

\bibitem{GuHo:83}
J.~Guckenheimer and P.~Holmes.
\newblock {\em Nonlinear {O}scillations, {D}ynamical {S}ystems and
  {B}ifurcations of {V}ector {F}ields}.
\newblock Springer-Verlag, New York, 1983.

\bibitem{IHS:98}
M.~Ipsen, F.~Hynne, and P.~G. S{\o}rensen.
\newblock Systematic derivation of amplitude equations and normal forms for
  dynamical systems.
\newblock {\em Chaos}, 8:834--852, 1998.

\bibitem{IHS:00}
M.~Ipsen, F.~Hynne, and {\mbox{P.G}.}~S{\o}rensen.
\newblock Amplitude equations for reaction-diffusion systems with a hopf
  bifurcation and slow real modes.
\newblock {\em Phys. D}, 136:66--92, 2000.

\bibitem{JeDe:86}
{\mbox{A.D}}.~Jepson and {\mbox{D.W}}.~Decker.
\newblock Convergence cones near bifurcation.
\newblock {\em SIAM J. Numer. Anal.}, 23:959--975, 1986.

\bibitem{KuGoDoDh:2005}
Yu.~A. Kuznetsov, W.~Govaerts, E.~J. Doedel, and A.~Dhooge.
\newblock Numerical periodic normalization for codim 1 bifurcations of limit
  cycles.
\newblock {\em SIAM J. Numer. Anal.}, 43:1407--1435, 2005.

\bibitem{Ku:99}
{\mbox{Yu.A}}.~Kuznetsov.
\newblock Numerical normalization techniques for all codim 2 bifurcations of
  equilibria in {ODEs}.
\newblock {\em SIAM J. Numer. Anal.}, 36:1104--1124, 1999.

\bibitem{Ku:04}
{\mbox{Yu.A}}.~Kuznetsov.
\newblock {\em Elements of {A}pplied {B}ifurcation {T}heory}.
\newblock Springer {V}erlag, Berlin, 2004.
\newblock Third Edition.

\bibitem{CONTENT}
{\mbox{Yu.A}}.~Kuznetsov and {\mbox{V.V}}.~Levitin.
\newblock {\sc content:} {A} multiplatform environment for analyzing dynamical
  systems.
\newblock ({\tt ftp.cwi.nl/pub/CONTENT}), 1995--1997.

\bibitem{KuMeVe:04}
{\mbox{Yu.A}}.~Kuznetsov, {\mbox{H.G.E}}.~Meijer, and L.~van Veen.
\newblock The fold-flip bifurcation.
\newblock {\em Int. J. Bif. Chaos}, 14:2253--2282, 2004.

\bibitem{Me:06}
{\mbox{H.G.E}.}~Meijer.
\newblock {\em Codimension 2 {B}ifurcations of {I}terated {M}aps}.
\newblock PhD thesis, Utrecht University, Netherlands, 2006.

\bibitem{ShNiNi:95}
A.~Shil'nikov, G.~Nicolis, and C.~Nicolis.
\newblock Bifurcation and predictability analysis of a low-order atmospheric
  circulation model.
\newblock {\em Int. J. Bif. Chaos}, 5:1701--1711, 1995.

\bibitem{Ve:03}
Lennaert van Veen.
\newblock Baroclinic flow and the {L}orenz-84 model.
\newblock {\em Int. J. Bif. Chaos}, 13:2117--2139, 2003.

\bibitem{WiCh:06}
S.~Wieczorek and {\mbox{W.W}}.~Chow.
\newblock Self-induced chaos in a single-mode inversionless laser.
\newblock {\em Phys. Rev. Lett.}, 97:113903, 2006.

\end{thebibliography}

\end{document}